\documentclass[12pt]{amsart}
\usepackage{amssymb,latexsym, amsmath, amscd, array, graphicx}

\swapnumbers
\numberwithin{equation}{section}

\theoremstyle{plain}

\newtheorem{thm}{Theorem}[section]
\newtheorem{theorem}[thm]{Theorem}

\newtheorem{prop}[thm]{Proposition}

\theoremstyle{definition}

\newtheorem{definition}[thm]{Definition}

\newtheorem{remark[thm]}{Remark}

\newtheorem{cor}[thm]{Corollary}

\newtheorem{quest}[thm]{Question}

\DeclareMathOperator{\cat}{{\mbox{\rm cat$_{\rm LS}$}}}

\def\Ord{\protect\operatorname{Ord}}

\def\scr{\mathcal}

\def\C{{\mathbb C}}

\def\1{\hbox{\rm\rlap {1}\hskip.03in{\rom I}}}
\def\Bbbone{{\rm1\mathchoice{\kern-0.25em}{\kern-0.25em}
{\kern-0.2em}{\kern-0.2em}I}}

\long\def\forget#1\forgotten{} %

\newcommand\ver[1]{\marginpar{\tiny Changed in Ver \VER}}

\begin{document}

\title {On the Lusternik-Schnirelmann category of spaces with\\
2-dimensional fundamental group}

\author[A.~Dranishnikov]{Alexander N. Dranishnikov}

\address{Department of Mathematics, University
of Florida, 358 Little Hall, Gainesville, FL 32611-8105, USA}
\email{dranish@math.ufl.edu}

\begin{abstract}
The following inequality
$$\cat X\le \cat Y+\left\lceil\frac{hd(X)-r}{r+1}\right\rceil$$
holds for every locally trivial fibration between $ANE$ spaces $f:X\to Y$
which admits a section and has the $r$-connected fiber where $hd(X)$ is the
homotopical dimension of $X$. We apply this inequality
to prove that
$$
\cat X\le \left\lceil\frac{\dim X-1}{2}\right\rceil+cd(\pi_1(X))$$ for every complex
$X$ with $cd(\pi_1(X))\le 2$.
\end{abstract}

\maketitle

\section{Introduction}
In \cite{DKR} we proved that if the Lusternik-Schnirelmann category
of a closed $n$-manifold, $n\ge 3$, equals 2 then the fundamental
group of $M$ is free. In opposite direction we proved that if the
fundamental group of  an $n$-manifold is free, then $\cat M\le n-2$.
Then J. Strom proved that $\cat X\le \frac{2}{3}n$ for every $n$-complex,
$n>4$, with free fundamental group \cite{St}. Yu. Rudyak suggested
that the coefficient $2/3$ in Strom's result should be improved to
$1/2$. Precisely, he conjectured that the function $f$ defined as
$f(n)=\max\{\cat M^n\}$ is asymptotically
$\frac{1}{2}n$ where the maximum is taken over all closed
$n$-manifolds with free fundamental group.

In this paper we prove Rudyak's conjecture. Our method gives the
same  estimate for $n$-complexes. Moreover, it gives the same
asymptotic upper bound for $\cat$ of $n$-complexes with the
fundamental group of cohomological dimension $\le 2$. In view of
this the following generalization of Rudyak's conjecture seems to be
natural.

\

CONJECTURE. {\em For every $k$ the function $f_k$ defined as
$$f_k(n)=\max\{\cat M^n\mid cd(\pi_1(M^n)\le k\}$$
is asymptotically $\frac{1}{2}n$.}

The smallest $k$ when it is unknown is 3.

The paper is organized as follows. Section 2 is an introduction to
the Lusternik-Schnirelmann category based on an analogy with the
dimension theory.  Section 3 contains a version of a fibration
theorem for $\cat$. In Section 4 this fibration theorem is applied
for the proof of Rudyak's conjecture.

\section{Kolmogorov-Ostrand's approach to the Lusternik-Schnirelmann category}
An open cover $\scr U=\{U_i\}$ of a topological space $X$ is
called {\em $X$-contractible} if each $U_i$ can be contracted to a
point in $X$. By the definition $\cat X\le n$ if there is an
$X$-contractible cover of $X$ by $n+1$ elements.

We recall \cite{CLOT} that a sequence $\emptyset=O_0\subset
O_1\subset\dots\subset O_{n+1}=X$ is called {\em categorical of
length $n+1$} if each difference $O_{i+1}\setminus O_i$ is
contained in an open set that contracts to a point. It was proven
in \cite{CLOT} that $\cat X\le n$ if and only if $X$ admits a
categorical sequence of length $n+1$.

Let $\scr U$ be a family of open set in a topological space $X$. The
multiplicity of $\scr U$ (or the order) at a point $x\in X$,
$\Ord_x\scr U$ is the number of elements of $\scr U$ that contain
$x$. The multiplicity of $\scr U$ is defined as $\Ord\scr
U=\sup_{x\in X}\Ord_x\scr U$. We recall that {\em the covering
dimension of a topological space $X$ does not exceed $n$, $\dim X\le
n$, if for every open cover $\scr C$ of $X$ there is an open
refinement $\scr U$ with $\Ord\scr U\le n+1$.}

The following proposition makes the LS-category analogous to
the covering dimension.
\begin{prop} For a paracompact topological space $X$,
$\cat X\le n$ if and only if $X$ admits an $X$-contractible locally
finite open cover $\scr V$ with $Ord\scr V\le n+1$.
\end{prop}
\begin{proof}
If $\cat X\le n$ then by the definition $X$ admits an open
contractible cover that consists of $n+1$ elements and therefore
its multiplicity is at most $n+1$.

Let $\scr V$ be a contractible cover of $X$ of multiplicity $\le
n+1$. We construct a categorical sequence $O_0\subset
O_1\subset\dots\subset O_{n+1}$ of length $n+1$. We define
$O_1=\{x\in X\mid \Ord_x\scr V\}=n+1$. Note that
$O_1=\cup_{V_i\in\scr V} V_0\cap\dots\cap V_n$. Note that this is
the disjoint union and every nonempty summand is $X$-contractible.
Thus $O_1$ is $X$-contractible. Next, we define $O_2=\{x\in X\mid
\Ord_x\scr V\ge n\}$. Then $O_2\setminus O_1=cup_{V_i\in\scr
V}(V_1\cap\dots\cap V_n\setminus O_1)$ is a disjoint union of closed
in $O_2$ subsets. Since this family is locally finite, we can take
open (in $O_2$ and hence in $X$) disjoint neighborhoods of these
summands $V_1\cap\dots\cap V_n\setminus O_1$ such that each of which
is contained in an open contractible set from $\scr V$. Define
$O_3=\{x\in X\mid Ord_x\ge n-1\}$ as the union of $n-1$-fold
intersections and so on. $O_{n+1}$ is the union of elements of $\scr
V$ (1-fold intersections) and hence $O_{n+1}=X$. Clearly, for every
$k$, $O_{k+1}\setminus O_k$ is the union of disjoint sets each of
which is contained in an element of $\scr V$. Thus, the categorical
sequence conditions are satisfied.
\end{proof}

A family $\scr U$ of subsets of $X$ is called a $k$-cover, $k\in
N$ if every subfamily of $k$ elements form a cover of $X$. The
following is obvious.
\begin{prop}\label{n-cover}
A family $\scr U$ of $m$-elements is an $(n+1)$-cover of $X$ if
and only if $\Ord_x\scr U\ge m-n$ for all $x\in X$.
\end{prop}
\begin{proof}
If $\Ord_x\scr U< m-n$ for some $x\in X$, then $n+1=m-(m-n)+1$
elements of $\scr U$ do not cover $x$.

If $n+1$ elements of $\scr U$ do not cover some $x$, then
$\Ord_x\scr U\le m-(n+1)<m-n$.
\end{proof}

Inspired by the work of Kolmogorov on Hilbert's 13 problem Ostrand
gave the following characterization of the covering dimension
\cite{Os}.
\begin{theorem}[Ostrand]
A metric space $X$ is of dimension $\leq n$ if and only if for
each open cover $\scr C$ of $X$ and each integer $m\geq n+1$,
there exist $m$ discrete families of open sets $\scr U_1,\cdots,\scr U_k$
such that their union $\cup\scr U_i$ is an  $(n+1)$-cover of $X$.
\end{theorem}

Let $\scr U$ be a family of subsets in $X$ and let $A\subset X$.
We denote by $\scr U|_A-\{U\cap A\mid U\in\scr U\}$.
\begin{definition} Let $f:X\to Y$ be a map.
An open cover $\scr U=\{U_0, U_1,\dots, U_n\}$ of $X$ is called {\em
uniformly $f$-contractible} if it satisfies the property that for
every $y\in Y$ there is a neighborhood $V$ such that the restriction
$\scr U|_{f^{-1}(V)}$ of $\scr U$ to the preimage $f^{-1}(V)$
consists of $X$-contractible sets.
\end{definition}

\begin{theorem}\label{criterion}
Let $\{U_0',\dots,U_n'\}$ be an open cover of a normal topological space $X$.
Then there is
an infinite open $(n+1)$-cover of $X$,  $\{U_k\}_{k=0}^{\infty}$
such that $U_k=U_k'$ for $k\le n$ and
$U_k=\cup_{i=0}^nV_i$ is the disjoint
union with $V_i\subset U_i$ for $k>n$.

In particular, if $\{U_0',\dots,U_n'\}$ is $X$-contractible, the
cover $\{U_k\}_{k=0}^{\infty}$ is $X$-contractible.
If $\{U_0',\dots,U_n'\}$ is uniformly $f$-contractible for some $f:X\to Y$, the
cover $\{U_k\}_{k=0}^{\infty}$ is uniformly $f$-contractible.
\end{theorem}
\begin{proof}
By induction on $m$ we construct the family
$\{U_i\}_{i=0}^{m}$. For $m=n$ we take $U_i=U_i'$.

Let $\scr U=\{ U_0,\dots, U_m\}$
be the corresponding family for $m\ge n$. Consider $Y=\{x\in X\mid \Ord_x\scr
U\le m-n\}$. In view of Proposition~\ref{n-cover} and the assumption, it
is a closed subset. We show that for every $i\le n$,  the set
$Y\cap U_i$ is closed in $X$. Let $x$ be a limit point of
$Y\cap U_i$ that does not belong to $U_i$. Let
$U_{i_0},\dots,U_{i_{m-n}}$ be the elements that contain $x$. Then
$\Ord_y\scr U=m-n+1$ for all $y\in Y\cap U_i\cap U_{i_0}\cap\dots\cap
U_{i_{m-n}}$. Contradiction.

We define recursively $F_0=Y\cap U_1$ and $F_{i+1}=Y\cap U_{i+1}\setminus
(\cup_{k=0}^iU_i)$. Note that $\{F_i\}$ is a disjoint finite family
of closed subsets.
We fix disjoint open neighborhoods $V_i$ of
$F_i$ with $V_i\subset U_i$. We define $U_{m+1}=\cup_iV_i$.
In view of
Proposition~\ref{n-cover}, $U_1,\dots, U_m,U_{m+1}$ is an $(n+1)$-cover.

Clearly, if all $U_i$ are $X$-contractible, $i\le n$,
then $U_{m+1}$ is $X$-contractible. If all $U_i$ are uniformly $f$-contractible,
for some $f:X\to Y$, then $U_{m+1}$ is uniformly $f$-contractible.
\end{proof}

\begin{cor}
For a normal topological space $X$, $\cat X\le n$ if and only if for any $m>n$
$X$ admits an open $(n+1)$-cover by $X$-contractible sets.
\end{cor}
This corollary is a $\cat$-analog of Ostrand's theorem. It also can
be found in \cite{CLOT} with further reference to \cite{Cu}.

\section{Fibration theorems for $\cat$}

\begin{definition}
The $\ast$-category $\cat^*f$ of a map $f:X\to Y$ is the minimal $n$,
if exists, such that there is a uniformly $f$-contractible
open cover $\scr U=\{U_0,
U_1,\dots, U_n\}$ of $X$.
\end{definition}

Note that $\cat^*c=\cat X$ for a constant map $c:X\to pt$. More
generally, $\cat^*\pi=\cat X$ for the projection $\pi:X\times Y\to
Y$.

\begin{theorem}\label{hur}
The inequality $\cat X\le\dim Y+\cat^*f$ holds true for any
continuous map.
\end{theorem}
\begin{proof}
The requirements to the spaces in the theorem are that the Ostrand
theorem holds true for $Y$, i.e. are fairly general (say, $Y$ is
normal).

Let $\dim Y=n$ and $\cat^*f=m$.  Let $\scr U=\{U_0,\dots U_m\}$ be
a uniformly $f$-contractible cover. For all $y\in Y$ denote by
$V_y$ a neighborhood of $y$ from the definition of the uniform
$f$-contractibility. In view of Ostrand's characterization of
dimension there is a refinement $\scr V=\scr V_0\cup\dots\cup\scr
V_{n+m}$ of the cover $\{V_y\mid y\in Y\}$ such that each family
$\scr V_i$ is disjoint. Let $V_i=\cup\scr V_i$. We apply Theorem
\ref{criterion} to extend the family $\scr U$ to a family
$\{U_0,\dots,U_{n+m}\}$. Consider the family $\scr
W=\{f^{-1}(V_i)\cap U_i\}$. Note that it is contractible in $X$.
We show that it is a cover of $X$. According to the Ostrand
theorem every $y\in Y$ is covered by $m+1$ elements
$V_{i_0},\dots, V_{i_m}$. By Theorem \ref{criterion} the family
$U_{i_0},\dots, U_{i_m}$ covers the fiber $f^{-1}(y)$.
\end{proof}
\begin{cor} [Corollary 9.35 \cite{CLOT}, \cite{OW}]
Let $p:X\to Y$ be a closed  map of $ANE$. If each fiber $p^{-1}(y)$
is contractible in $X$, then $\cat X\le\dim Y$.
\end{cor}
\begin{proof} In this case $\cat^*p=0$. Indeed,
since $X$ is an ANE, a contraction of $p^{-1}(y)$ to a point can be extended
to a neighborhood $U$. Since the map $p$ is closed there is a neighborhood
$V$ of $y$ such that $p^{-1}(V)\subset U$.\qed
\end{proof}

We recall that the {\em homotopical dimension} of a space $X$, $hd(X)$, is the
minimal dimension of a $CW$-complex homotopy equivalent to $X$ \cite{CLOT}.

The elementary obstruction theory implies the following.
\begin{prop}\label{section}
Let $p:E\to X$ be a fibration with $n-1$-connected fiber where
$n=hd(X)$. Then $p$ admits a section.
\end{prop}

We recall that $\lceil x\rceil$ denotes the smallest integer $n$
such that $x\le n$.
\begin{theorem}\label{r-connectedhur}
Suppose that a locally trivial fibration $f:X\to Y$ with an
$r$-connected fiber $F$ admits a section. Then
$$\cat^*f\le\left\lceil\frac{hd(X)-r}{r+1}\right\rceil.$$
Moreover,
$$\cat X\le \cat Y+
\left\lceil\frac{hd(X)-r}{r+1}\right\rceil.$$
\end{theorem}
\begin{proof}
Let $\cat Y=m$ and  $hd(X)=n$.

Let $s:Y\to X$ be a section. We apply the fiber-wise Serre
construction to $f$ with this $s(y)$ as the base points. Then we
apply the fiber-wise Ganea construction to obtain a map
$\xi_k:E_k\to X$. Namely,
$$
E_0=\{\phi\in X^I\mid |f(\phi(I))|=1,\ \phi(0)=sf\phi(I)\}
$$
with $\xi_0(\phi)=\phi(1)$ and
$$
E_{k+1}=\{\phi\ast\psi\in E_0\ast E_k\mid
\xi_0(\phi)=\xi_k(\psi)\}
$$
with the natural projection $\xi_{k+1}:E_{k+1}\to X$. Note that
$\xi_k$ is a Hurewicz fibration with the fiber the join product
$\ast^{k+1}\Omega F$ of
$k+1$ copies of the loop space $\Omega F$. Thus, it is
$(k+(k+1)r-1)$-connected. By Proposition \ref{section} there is a section
$\sigma:X\to E_k$ whenever $k(r+1)+r\ge n$. The smallest such $k$ is
equal to $\lceil\frac{gcd(X)-r}{r+1}\rceil$.

For each $x\in X$ the element$\sigma(x)$ of $\ast_k\Omega F$ can be
presented as the $(k+1)$-tuple
$$\sigma(x)=((\phi_0,t_0),\dots,(\phi_k,t_k))\mid \sum t_i=1, t_i\ge
0).$$ We use the notation $\sigma(x)_i=t_i$. Clearly, $\sigma(x)_i$
is a continuous function.

A section $\sigma:X\to E_k$ defines a cover $\scr U=\{U_0,\dots,
U_k\}$ of $X$ as follows:
$$
U_i=\{x\in X\mid\sigma(x)_i>0\}.
$$
Let $\{U_0,\dots, U_{m+k}\}$ be an extension of $\scr U$ to an
$(k+1)$-cover of $X$ from Theorem \ref{criterion}.

Let $\scr V=\{V_0,\dots V_{m+k}\}$ be an open $Y$-contractible $(m+1)$-cover
of $Y$. We show that the sets $W_i=f^{-1}(V_i)\cap U_i$ are
contractible in $X$ for all $i$.
By the construction of $U_i$ for $i\le n$ for every $x\in U_i$
there is a canonical path connecting $x$ with $sf(x)$. We use
these paths to contract $f^{-1}(V_i)\cap U_i$ to $s(V_i)$ in $X$.
Then we apply a contraction of $s(V_i)$ to a point in $s(Y)$. If
$i>k$, we apply the last condition of the first part of Theorem \ref{criterion}.

Similarly as in the proof of Theorem \ref{hur} we show that
$\{W_i\}_{i=0}^{m+k}$ is a cover of $X$. Since $\scr V$ is an
$(m+1)$-cover, by Proposition \ref{n-cover} every $y\in Y$ is
covered by at least $k+1$ elements $V_{i_0},\dots, V_{i_k}$ of $\scr
V$. By the construction $U_{i_0},\dots, U_{i_k}$ is a cover of $X$.
Hence $W_{i_0},\dots, W_{i_k}$ covers $f^{-1}(y)$.
\end{proof}

\section{The Lusternik-Schnirelmann category of complexes with low dimensional fundamental groups}

\begin{theorem}
For every complex $X$ with $cd(\pi_1(X))\le 2$ the following inequality
holds true:
$$
\cat X\le cd(\pi_1(X)) +\left\lceil\frac{hd(X)-1}{2}\right\rceil.$$
\end{theorem}
\begin{proof}
Let $\pi=\pi_1(X)$.
We consider Borel's construction
$$
\begin{CD}
\tilde X @<<< \tilde X\times E\pi @>>> E\pi\\
@VVV @VVV @VVV\\
X @<g<< \tilde X\times_{\pi}E\pi @>f>> B\pi.\\
\end{CD}
$$
We claim that there is a section $s:B\pi\to \tilde X\times_{\pi}E\pi$ of $f$.
By the condition $cd\pi\le 2$ we may assume that $B\pi$ is a complex
of dimension $\le 3$. Since a fiber of $f$ is simply connected,
there is a lift of the 2-skeleton. The condition $cd\pi\le 2$
implies $H^3(B\pi,E)=0$ for every
$\pi$-module. Thus, we have no obstruction for the lift of the 3-skeleton.

We apply Theorem \ref{r-connectedhur} to obtain the inequality
$$
\cat X\le \cat(B\pi) +\left\lceil\frac{hd(\tilde X\times_{\pi}E\pi)-1}{2}\right\rceil.$$
Since $g$ is a fibration with the homotopy trivial fiber, the space
$\tilde X\times_{\pi}E\pi$ is
homotopy equivalent to $X$.
Thus, $hd(\tilde X\times_{\pi}E\pi)=hd(X)$. Note that $\cat B\pi=cd\pi$.
\end{proof}

\begin{cor}
For every complex $X$ with free fundamental group:
$$
\cat X\le 1 +\left\lceil\frac{\dim X-1}{2}\right\rceil.$$
\end{cor}
Note that this estimate is sharp on $X=S^1\times \C P^n$.
\begin{cor}
For every 3-dimensional complex $X$ with free fundamental group:
$\cat X\le 2.$
\end{cor}
This Corollary can be also derived from the fact that in the case
of free fundamental group every 2-complex is homotopy equivalent
to the wedge of circles and 2-spheres \cite{KR}.

It is unclear whether the estimate $ \cat X\le 2+\lceil\frac{\dim
X-1}{2}\rceil$ is sharp for complexes with $cd(\pi_1(X))=2$. It is
sharp when the answer to the following question is affirmative.

\begin{quest} Does there exists a 4-complex $K$ with free fundamental group
and with $\cat (K\times S^1)=4$?
\end{quest}

Indeed, for $X=K\times S^1$ we would have the equality
$4=2+\lceil\frac{5-1}{2}\rceil$. Note that $cd(\pi_(X))=2$. There is
a connection between this question and the Problem 1.5 from
\cite{DKR} which asks whether $\cat M^4\le 2$ for closed
$4$-manifolds with free fundamental group.

\end{document}